\documentclass[11pt]{article}
\usepackage{amsfonts,amssymb,mathrsfs, latexsym}
\usepackage{amsmath,amscd, eufrak}
\usepackage{color}
\usepackage[all]{xy}
\definecolor{red}{rgb}{1,0,0}
\definecolor{blue}{rgb}{0,0,1}
\definecolor{darkblue}{rgb}{0,0,0.4}

\begin{document}
\newcommand{\dickebox}{{\vrule height5pt width5pt depth0pt}}
\newtheorem{Def}{Definition}[section]
\newtheorem{Bsp}{Example}[section]
\newtheorem{Prop}[Def]{Proposition}
\newtheorem{Theo}[Def]{Theorem}
\newtheorem{Lem}[Def]{Lemma}
\newtheorem{Koro}[Def]{Corollary}
 \newcommand{\lra}{\longrightarrow}
 \newcommand{\ra}{\rightarrow}
 \newcommand{\add}{{\rm add\, }}
\newcommand{\gd}{{\rm gl.dim\, }}
\newcommand{\End}{{\rm End\, }}
\newcommand{\overpr}{$\hfill\square$}
\newcommand{\rad}{{\rm rad\,}}
\newcommand{\soc}{{\rm soc\,}}
\renewcommand{\top}{{\rm top\,}}
\newcommand{\pd}{{\rm proj.dim\, }}

\newcommand{\cpx}[1]{#1^{\bullet}}
\newcommand{\D}[1]{{\mathscr D}(#1)}
\newcommand{\Dz}[1]{{\rm D}^+(#1)}
\newcommand{\Df}[1]{{\rm D}^-(#1)}
\newcommand{\Db}[1]{{\mathscr D}^b(#1)}
\newcommand{\C}[1]{{\mathscr C}(#1)}
\newcommand{\Cz}[1]{{\rm C}^+(#1)}
\newcommand{\Cf}[1]{{\rm C}^-(#1)}
\newcommand{\Cb}[1]{{\mathscr C}^b(#1)}
\newcommand{\K}[1]{{\mathscr K}(#1)}
\newcommand{\Kz}[1]{{\rm K}^+(#1)}
\newcommand{\Kf}[1]{{\rm K}^-(#1)}
\newcommand{\Kb}[1]{{\mathscr K}^b(#1)}
\newcommand{\modcat}[1]{#1\mbox{{\rm -mod}}}
\newcommand{\stmodcat}[1]{#1\mbox{{\rm -{\underline{mod}}}}}
\newcommand{\pmodcat}[1]{#1\mbox{{\rm -proj}}}
\newcommand{\imodcat}[1]{#1\mbox{{\rm -inj}}}
\newcommand{\opp}{^{\rm op}}
\newcommand{\otimesL}{\otimes^{\rm\bf L}}
\newcommand{\rHom}{{\rm\bf R}{\rm Hom}\,}
\newcommand{\projdim}{\pd}
\newcommand{\Hom}{{\rm Hom \, }}
\newcommand{\Coker}{{\rm coker}\,\,}
\newcommand{\Ext}{{\rm Ext}}
\newcommand{\StHom}{{\rm \underline{Hom} \, }}
{\Large \bf
\begin{center}
Almost $\cal D$-split sequences and derived equivalences
\end{center}}
\medskip

\centerline{\bf Wei Hu and Changchang Xi$^*$}
\begin{center} School of Mathematical Sciences, Beijing Normal University, \\
Laboratory of Mathematics and Complex Systems, MOE, \\
100875 Beijing, People's Republic of  China \\ E-mail: xicc@bnu.edu.cn \quad hwxbest@163.com\\
\end{center}

\renewcommand{\thefootnote}{\alph{footnote}}
\setcounter{footnote}{-1} \footnote{ $^*$ Corresponding author.
Email: xicc@bnu.edu.cn; Fax: 0086 10 58802136; Tel.: 0086 10
58808877.}
\renewcommand{\thefootnote}{\alph{footnote}}
\setcounter{footnote}{-1} \footnote{2000 Mathematics Subject
Classification: 16G70,18E30;16G10,18G20.}
\renewcommand{\thefootnote}{\alph{footnote}}
\setcounter{footnote}{-1} \footnote{Keywords: almost $\cal D$-split
sequence, Auslander-Reiten triangle, BB-tilting module, derived
equivalence, stable \newline equivalence.}

\abstract{In this paper, we introduce almost $\cal D$-split
sequences and establish an elementary but somewhat surprising
connection between derived equivalences and Auslander-Reiten
sequences via BB-tilting modules. In particular, we obtain derived
equivalences from Auslander-Reiten sequences (or $n$-almost split
sequences), and Auslander-Reiten triangles.}

\section{Introduction}
Derived equivalence and Auslander-Reiten sequence are two important
objects in the modern representation theory of algebras and groups.
On the one hand, derived equivalence preserves many significant
invariants of groups and algebras; for example, the number of
irreducible representations, Cartan determinants, Hochschild
cohomology groups, algebraic K-theory and G-theory(see
\cite{Broue1994}, \cite{HappelTri} and \cite{DuggerShipley}). One of
the fundamental results on derived categories may be the Morita
theory for derived categories established by Rickard in his several
papers \cite{RickMoritaTh, RickDstable, RickDFun}, which says that
two rings $A$ and $B$ are derived-equivalent if and only if there is
a tilting complex $T$ of $A$-modules such that $B$ is isomorphic to
the endomorphism ring of $T$. Thus, starting with a ring $A$, we may
construct theoretically all rings which are derived-equivalent to
$A$ by finding all tilting complexes of $A$-modules. However, in
practice, it is not easy to show that two given rings are
derived-equivalent by finding a suitable tilting complex, as is
indicated by the famous unsolved Broue's abelian defect group
conjecture, which states that the module categories of a block
algebra $A$ of a finite group algebra and its Brauer correspondent
$B$ should have equivalent derived categories if their defect groups
are abelian (see \cite{Broue1994}). On the other hand, as is
well-known, Auslander-Reiten sequence is of significant importance
in the modern representation theory of Artin algebras, it contains
rich combinatorial information on the module category (see
\cite{AusReiten}). A natural and fundamental question is: Is there
any relationship between Auslander-Reiten sequences and derived
equivalences ? In other words, is it possible to construct derived
equivalences from Auslander-Reiten sequences or $n$-almost split
sequences or Auslander-Reiten triangles ?

In the present paper, we shall provide an affirmative answer to this
question  and construct derived equivalences by the so-called almost
$\cal D$-split sequences (see Definition \ref{def1} below). Such
sequences include Auslander-Reiten sequences and occur very
frequently in the representation theory of Artin algebras (see the
examples in Section \ref{sect3} below). Our result in this direction
can be stated in the following general form:

\begin{Theo}
Let $\cal C$ be an additive category and $M$ be an object in $\cal
C$. Suppose
$$X\longrightarrow M'\longrightarrow Y
$$ is an almost \emph{add}$(M)$-split sequence in
$\cal C$. Then the endomorphism ring $\End_{\cal C}(M\oplus X)$ of
$M\oplus X$ and the endomorphism ring $\End_{\cal C}(M\oplus Y)$ of
$M\oplus Y $ are derived-equivalent via a tilting module. Moreover,
the finitistic dimension of $\End_{\cal C}(M\oplus X)$ is finite if
and only if so is the finitistic dimension of $\End_{\cal C}(M\oplus
Y)$. \label{thm1}
\end{Theo}

This result reveals a mysterious connection between Auslander-Reiten
sequences and derived equivalences, namely we have the following
corollary.

\begin{Koro}
Let $A$ be an Artin algebra.

$(1)$ Suppose $0\lra X_i\lra M_i\lra X_{i-1}\lra 0$ is an
Auslander-Reiten sequence of finitely generated $A$-modules for
$i=1,2,\cdots, n$. Let $M=\bigoplus_{i=1}^nM_i$. Then
$\End_A(M\oplus X_n)$ and $\End_A(M\oplus X_0)$ are
derived-equivalent via an $n$-BB-tilting module. In particular, if
$\; 0\lra X\lra M\longrightarrow Y\lra 0$ is an Auslander-Reiten
sequence, then the endomorphism algebras $\End_A(X\oplus M)$ and
$\End_A(M\oplus Y)$ are derived-equivalent via BB-tilting module,
and have the same Cartan determinant.

$(2)$ If $A$ is self-injective and $X$ is an $A$-module, then the
endomorphism algebra $\End(A\oplus X)$ of $A\oplus X$ and the
endomorphism algebra $\End_A(A\oplus \Omega(X))$ of $A\oplus
\Omega(X)$ are derived-equivalent, where $\Omega$ is the syzygy
operator. \label{cor1}
\end{Koro}

Thus, by Corollary \ref{cor1} or more generally, by Proposition
\ref{n-ars}  in Section \ref{sect3} below, one can produce a lot of
derived equivalences from Auslander-Reiten sequences or $n$-almost
split sequences. We stress that the algebra $\End_A(X\oplus M)$ and
the algebra $\End_A(M\oplus Y)$ in Corollary \ref{cor1} may be very
different from each other (see the examples in Section \ref{sect7}),
though the mesh diagram of the Auslander-Reiten sequence is somehow
symmetric. Another result related to Corollary \ref{cor1} is
Proposition \ref{artriangle} in Section \ref{sectriangle} below,
which produces derived equivalences from Auslander-Reiten triangles
in a triangulated category. In particular, we have

\begin{Koro}   Let $A$ be a self-injective
Artin algebra. Suppose $0\lra X\lra M\lra Y\lra 0$ is an
Auslander-Reiten sequence such that
$\Omega^{-1}(X)\not\in\add(M\oplus Y)$. Then
$\underline{\End}_A(M\oplus X)$ and $\underline{\End}_A(M\oplus Y)$
are derived-equivalent, where $\underline{\End}_A(M)$ denotes the
stable endomorphism algebra of an $A$-module $_AM$. \label{stable}
\end{Koro}

\medskip
The paper is organized as follows: In Section \ref{sect2}, we recall
briefly some basic notions and a fundamental result of Rickard on
derived categories. Our main results, Theorem \ref{thm1}, is proved
in Section \ref{sect3}, where we also provide several
generalizations of Corollary \ref{cor1}; among others is a
formulation of Corollary \ref{cor1}(1) for $n$-almost split
sequences. In section \ref{sect3+}, we point out that if an almost
$\cal D$-split sequence is given by an Auslander-Reiten sequence
then Theorem \ref{thm1} can be viewed as a $``$generalized" version
of a BB-tilting module. Thus an $n$-almost split sequence or
concatenating $n$ Auslander-Reiten sequences provides us a natural
way to get an $n$-BB-tilting module (for definition, see Section
\ref{sect3+}). In Section \ref{sectriangle}, we discuss how to get
derived equivalences from Auslander-Reiten triangles in a
triangulated category. In particular, Corollary \ref{stable} is
proved in this section. In the last section we present an example to
illustrate our main result.

\section{Preliminaries \label{sect2}}
In this section, we recall some basic definitions and results
required in our proofs.

Let $\cal C$ be an additive category. For two morphisms $f:X\lra Y$
and $g:Y\lra Z$ in $\cal C$, the composition of $f$ with $g$ is
written as $fg$, which is a morphism from $X$ to $Z$. But for two
functors $F:\mathcal{C}\lra \mathcal{D}$ and
$G:\mathcal{D}\lra\mathcal{E}$ of categories, their composition is
denoted by $GF$. For an object $X$ in $\mathcal{C}$, we denote by
$\add(X)$ the full subcategory of $\cal C$ consisting of all direct
summands of finite sums of copies of $X$.

A complex $\cpx{X}$ over $\cal C$ is a sequence of morphisms
$d_X^{i}$ between objects $X^i$ in $\cal C$: $ \cdots \ra
X^{i-1}\stackrel{d_X^{i-1}}{\lra} X^i\stackrel{d_X^i}{\lra}
X^{i+1}\stackrel{d_X^{i+1}}{\lra}X^{i+2}\ra\cdots$, such that
$d_X^id_X^{i+1}=0$ for all $i \in {\mathbb Z}$. We write
$\cpx{X}=(X^i, d_X^i)$. The category of all complexes over $\cal C$
with the usual complex maps of degree zero is denoted by $\C{\cal
C}$. The homotopy and derived categories of complexes over
$\mathcal{C}$ are denoted by $\K{\mathcal C}$ and $\D{\mathcal{C}}$,
respectively. The full subcategory of $\C{\cal C}$ consisting of
bounded complexes over $\mathcal{C}$ is denoted by $\Cb{\mathcal
C}$. Similarly, $\Kb{\mathcal C}$ and $\Db{\mathcal{C}}$ denote the
full subcategories consisting of bounded complexes in $\K{\mathcal
C}$ and $\D{\mathcal C}$, respectively.

An object $X$ in a triangulated category $\mathcal{C}$ with a shift
functor $[1]$ is called \emph{self-orthogonal} if
$\Hom_{\mathcal{C}}(X,X[n])=0$ for all integers $n\neq 0$.

Let $A$ be a ring with identity. By $A$-module we shall mean a left
$A$-module. We denote by $A$-Mod the category of all $A$-modules, by
$A$-mod the category of all finitely presented $A$-modules, and by
$\pmodcat{A}$ (respectively, $A$-inj) the category of finitely
generated projective ( respectively, injective) $A$-modules. Let $X$
be an $A$-module. If $f: P\lra X$ is a projective cover of $X$ with
$P$ projective, then the kernel of $f$ is called a \emph{syzygy} of
$X$, denoted by $\Omega(X)$. Dually, if $g: X\lra I$ is an injective
envelope with $I$ injective, then the cokernel of $g$ is called a
\emph{co-syzygy} of $X$, denoted by $\Omega^{-1}(X)$. Note that a
syzygy or a co-syzygy of an $A$-module $X$ is determined, up to
isomorphism, uniquely by $X$. Hence we may speak of the syzygy and
the co-syzygy of a module.

It is well-known that $\K{A\mbox{{\rm -Mod}}}$, $\Kb{A\mbox{{\rm
-Mod}}}$, $\D{A\mbox{{\rm -Mod}}}$ and $\Db{A\mbox{{\rm -Mod}}}$ all
are triangulated categories. Moreover, it is known that if $X \in
\Kb{A\mbox{-proj}}$ or $Y \in \Kb{A\mbox{-inj}}$, then
$\Hom_{\Kb{A-\mbox{Mod}}}(X, Z)$ $\simeq \Hom_{\Db{A-\mbox{Mod}}}(X,
Z)$ and Hom$_{\Kb{A-\mbox{Mod}}}(Z, Y) \simeq$
Hom$_{\Db{A\mbox{-Mod}}}(Z, Y)$ for all $Z\in \Db{A\mbox{{\rm
-Mod}}}$.

For further information on triangulated categories, we refer to
\cite{HappelTri}. In \cite{RickMoritaTh}, Rickard proved the
following theorem.

\begin{Theo}
For two rings $A$ and $B$ with identity, the following are
equivalent:

$(a)$ $\Db{A\mbox{{\rm -Mod}}}$ and $\Db{B\mbox{{\rm -Mod}}}$ are
equivalent as triangulated categories;

$(b)$ $\Kb{\pmodcat{A}}$ and $\Kb{\pmodcat{B}}$ are equivalent as
triangulated categories;

$(c)$ $B\simeq \End_{\Kb{\pmodcat{A}}}(\cpx{T})$, where $\cpx{T}$ is
a complex in $\Kb{\pmodcat{A}}$ satisfying

\parindent=1cm $(1)$  $\cpx{T}$ is self-orthogonal in $\Kb{\pmodcat{A}}$,

$(2)$ $\add(\cpx{T})$ generates $\Kb{\pmodcat{A}}$ as a triangulated
category. \label{rickard}
\end{Theo}

If two rings $A$ and $B$ satisfy the equivalent conditions of
Theorem \ref{rickard}, then $A$ and $B$ are said to be
\emph{derived-equivalent}. A complex $\cpx{T}$ in $\Kb{\pmodcat{A}}$
satisfying the conditions $(1)$ and $(2)$ in Theorem \ref{rickard}
is called a \emph{tilting complex} over $A$. Given a derived
equivalence $F$ between $A$ and $B$, there is a unique (up to
isomorphism) tilting complex $\cpx{T}$ over $A$ such that
$F\cpx{T}=B$. This complex $\cpx{T}$ is called a tilting complex
\emph{associated} to $F$.

\medskip
To get derived equivalences, one may use tilting modules. Recall
that a module $T$ over a ring $A$ is called a \emph{tilting module}
if

$(1)$ $T$ has a finite projective resolution $0\lra
P_n\lra\cdots\lra P_0\lra T\lra 0$, where each $P_i$ is a finitely
generated projective $A$-module;

$(2)$ $\Ext^i_A(T,T)=0$ for all $i>0$, and

$(3)$ there is an exact sequence $0\lra A\lra T^0\lra\cdots\lra
T^m\lra 0$ of $A$-modules with each $T^i$ in $\add(T)$.

\smallskip
It is well-known that each tilting module supplies a derived
equivalence. The following result in \cite{CPS} is a generalization
of a result in \cite[Theorem 2.10]{HappelTri}.

\begin{Lem} Let $A$ be a ring, $_AT$ a tilting $A$-module and $B=\End_A(T)$. Then $A$ and
$B$ are derived-equivalent. In this case, we say that $A$ and $B$
are derived-equivalent via a tilting module.
\label{tiltingmoduleforrings}
\end{Lem}

In Theorem \ref{rickard}, if both $A$ and $B$ are left coherent
rings, that is, rings for which the kernels of any homomorphisms
between finitely generated projective modules are finitely
generated, then $A\mbox{{\rm -mod}}$ and $B\mbox{{\rm -mod}}$ are
abelian categories, and the equivalent conditions in Theorem
\ref{rickard} are further equivalent to the condition

\medskip
{\it $(d)$ $\Db{A\mbox{{\rm -mod}}}$ and $\Db{B\mbox{{\rm -mod}}}$
are equivalent as triangulated categories.}

\medskip
A special class of coherent rings is the class of Artin algebras.
Recall that an \emph{Artin} $R$-\emph{algebra}  over a commutative
Artin ring $R$ is an $R$-algebra $A$ such that $A$ is a finitely
generated $R$-module.  For the module category over an Artin
algebra, there is the notion of Auslander-Reiten sequence, or
equivalently, almost split sequence. It plays an important role in
the modern representation theory of algebras and groups. Recall that
a short exact sequence $0\lra X\stackrel{f}{\lra}
Y\stackrel{g}{\lra} Z\lra 0$ in $A$-mod is called an
\emph{Auslander-Reiten sequence} if

(1) the sequence does not split,

(2) $X$ and $Z$ are indecomposable,

(3) for any morphism $h: V\lra Z$ in $A$-mod, which is not a split
epimorphism, there is a homomorphism $f':V\lra Y$ in $A$-mod such
that $h=f'f$, and

(4) for any morphism $h: X\lra V$ in $A$-mod, which is not a split
monomorphism, there is a homomorphism $f': Y\lra V$ in $A$-mod such
that $h=ff'$.

\medskip
For an introduction to Auslander-Reiten sequences and
representations of Artin algebras, we refer the reader to the
excellent book \cite{AusReiten}.

\section{Almost $\cal D$-split sequences and derived equivalences \label{sect3}}
In this section, we shall construct derived equivalences from
Auslander-Reiten sequences. This builds a linkage between
Auslander-Reiten sequences (or n-almost split sequences) and derived
equivalences. We start first with a general setting by introducing
the notion of almost $\cal D$-split sequences, which is a slight
generalization of Auslander-Reiten sequences, and then use these
sequences to construct derived equivalences between the endomorphism
rings of modules involved in almost $\cal D$-split sequences. In
Section \ref{sectriangle}, we shall consider the question of getting
derived equivalences from Auslander-Reiten triangles.

Now we recall some definitions from \cite{AS}.

Let $\cal C$ be a category, and let $\cal D$ be a full subcategory
of $\cal C$, and $X$ an object in $\cal C$. A morphism $f: D\lra X$
in $\cal C$ is called a \emph{right} $\cal D$-\emph{approximation}
of $X$ if $D\in {\cal D}$ and the induced map Hom$_{\cal C}(-,f)$:
Hom$_{\cal C}(D',D)\lra$ Hom$_{\cal C}(D',X)$ is surjective for
every object $D'\in {\cal D}$.  A morphism $f:X\lra Y$ in $\cal C$
is called \emph{right minimal} if any morphism $g: X\lra X$ with
$gf=f$ is an automorphism. A minimal right $\cal D$-approximation of
$X$ is a right $\cal D$-approximation of $X$, which is right
minimal. Dually, there is the notion of a \emph{left} $\cal
D$-\emph{approximation} and a \emph{minimal left} $\cal
D$-\emph{approximation}. The subcategory $\cal D$ is called
\emph{contravariantly} (respectively, \emph{covariantly} )
\emph{finite} in $\cal C$ if every object in $\cal C$ has a right
(respectively, left) $\cal D$-approximation. The subcategory $\cal
D$ is called \emph{functorially finite} in $\cal C$ if $\cal D$ is
both contravariantly and covariantly finite in $\cal C$.

Let $\cal C$ be an additive category and $e: X\lra X$  an idempotent
morphism in $\cal C$. We say that $e$ \emph{splits} if there are
objects $X'$ and $X''$ in $\cal C$ and an isomorphism $\varphi:
X'\oplus X''\lra X$ such that $\varphi e = \pi\lambda\varphi$, where
$\pi: X'\oplus X''\lra X'$ and $\lambda: X'\lra X'\oplus X''$ are
the canonical morphisms. In an arbitrary additive category, all
idempotents need not split, but of course, in the case where $\cal
C$ is an abelian category, every idempotent splits. If all
idempotents in $\cal C$ split, then so is every full subcategory
$\cal D$ of $\cal C$ which is closed under direct summands.
Moreover, for an additive category  $\cal C$ such that every
idempotent splits, we know that, for each object $M$ in $\cal C$,
the functor $\Hom_{\cal C}(M, -)$ induces an equivalence between
$\add(M)$ and $\End_{\cal C}(M)$-proj.

\begin{Def} Let $\cal C$ be  an additive category and $\cal D$ a full
subcategory of $\cal C$. A sequence
$$X\stackrel{f}{\longrightarrow}M\stackrel{g}{\longrightarrow}Y$$ in $\cal C$
is called an almost $\cal D$-split sequence if

$(1)$ $M\in {\cal D}$;

$(2)$ $f$ is a left $\cal D$-approximation of $X$, and $g$ is a
right $\cal D$-approximation of $Y$;

$(3)$ $f$ is a kernel of $g$, and $g$ is a cokernel of $f$.
\label{def1}
\end{Def}

Recall that a morphism $f: Y\lra X$ in an additive category $\cal C$
is a \emph{kernel} of a morphism $g:X\lra Z$ in $\cal C$ if $fg=0$,
and for any morphism $h: V\lra X$ in $\cal C$ with $hg=0$, there is
a unique morphism $h':V\lra Y$ such that $h=h'f$. Note that if a
morphism has a kernel in $\cal C$ then it is unique up to
isomorphism. A \emph{cokernel} of a given morphism in $\cal C$ is
defined dually. If $f: Y\lra X$ in $\cal C$ is a kernel of a
morphism $g:X\lra Z$ in $\cal C$, then $f$ is a monomorphism, that
is, if $h_i: U\lra Y$ is a morphism in $\cal C$ for $i=1,2,$ such
that $h_1f=h_2f$, then $h_1=h_2$. Similarly, if $g: X \lra Z$ in
$\cal C$ is a cokernel of a morphism $f: Y\lra X$ in $\cal C$, then
$g$ is an epimorphism, that is, if $h_i: Z\lra V$ is a morphism in
$\cal C$ for $i=1,2,$ such that $gh_1=gh_2$, then $h_1=h_2$.

Notice that an almost $\cal D$-split sequence may split, whereas an
Auslander-Reiten sequence never splits. Now we give some examples of
almost $\cal D$-split sequences.

\smallskip
{\bf Examples.} (a) Let $A$ be an Artin algebra and ${\cal C}$ =
$A$-mod. Suppose ${\cal D}$ is the full subcategory of $A$-mod
consisting of all projective-injective $A$-modules in $\cal C$. If
$g: M\lra X$ is a surjective homomorphism in $A$-mod with $M\in \cal
D$, then the sequence  $ 0\lra \mbox{ker}(g)\lra M\lra X\lra 0$ is
an almost $\cal D$-split sequence in $\cal C$, where ker$(g)$ stands
for the kernel of the homomorphism $g$.

(b)  Let $A$ be an Artin algebra and ${\cal C}$ = $A$-mod. Suppose
$0\lra X\lra M\lra Y\lra 0$ is an Auslander-Reiten sequence. Let $N$
be any module such that $M\in $ add($N$), but neither $X$ nor $Y$
belongs to add$(N)$. If we take $\cal D$ = add$(N)$, then the
Auslander-Reiten sequence is an almost $\cal D$-split sequence in
$\cal C$.

(c) Let $A$ be an Artin algebra and $M \in A\modcat$. Recall that
$M$ is an \emph{almost complete tilting module} if $M$ is a partial
tilting module (that is, $M$ has finite projective dimension and
$\Ext_A^i(M,M)=0$ for all $i>0$), and if the number of all
non-isomorphic direct summands of $M$ equals the number of
non-isomorphic simple $A$-modules minus $1$. An indecomposable
$A$-module $X \in A\modcat$ is called a \emph{tilting complement} to
$M$ if $M\oplus X$ is a tilting $A$-module. If an almost complete
tilting module $M$ is faithful, then there is an exact (not
necessarily infinite) sequence
$$ 0 \lra X_0\stackrel{f_1}{\lra} M_1\stackrel{f_2}{\lra} M_2 \stackrel{f_3}{\lra}
\cdots $$ of $A$-modules such that $M_i\in$ add$(M)$. Moreover, if
we define $X_i = \mbox{coker}(f_i)$, the co-kernel of $f_i$ for
$i\ge 1$, then $X_i\not\simeq X_j$ for $i \neq j$,
proj.dim$_A(X_i)\ge i$ for any $i$, and $\{X_i\mid i\ge 0\}$ is a
complete set of non-isomorphic indecomposable tilting complements to
$M$. In addition, each $X_i\lra M_{i+1}$ is a minimal left
add$(M)$-approximation of $X_i$ and each $M_j \lra X_j$ is a minimal
right add$(M)$-approximation of $X_j$. Thus the sequence $0\lra
X_i\lra M_{i+1}\lra X_{i+1}\lra 0$ is an almost add$(M)$-split
sequence in $A$-mod for all $i\ge 0$. For further information on
almost complete tilting modules and relationship with the
generalized Nakayama conjecture, we refer the reader to \cite{bs98}
and \cite{hu98}.

\medskip
Now we consider some properties of an almost $\cal D$-split
sequence.

\begin{Prop} Let $\cal C$ be  an additive category and $\cal D$ a full
subcategory of $\cal C$.

$(1)$ Suppose $\cal D'$ is a full subcategory of $\cal D$. If a
sequence $X\longrightarrow M\longrightarrow Y$ in $\cal C$ is an
almost $\cal D$-split sequence with $M\in {\cal D}'$, then it is an
almost ${\cal D}'$-split sequence in $\cal C$.

$(2)$ If $X\lra M\stackrel{g}{\lra}Y$ and $X'\lra
M'\stackrel{g'}{\lra} Y'$ are almost $\cal D$-split sequences in
$\cal C$ such that both $g$ and $g'$ are right minimal, then
$Y\simeq Y'$ if and only if the two sequences are isomorphic.
Similarly, If $X\stackrel{f}{\lra} M\lra Y$ and
$X'\stackrel{f'}{\lra} M'\lra Y'$ are almost $\cal D$-split
sequences in $\cal C$ such that both $f$ and $f'$ are left minimal,
then $X\simeq X'$ if and only if the two sequences are isomorphic.
\end{Prop}

{\it Proof.} (1) is clear. We prove the first statement of (2). If
the two sequences are isomorphic, then $X\simeq X'$ and $Y\simeq
Y'$. Now assume that $\phi: Y\lra Y'$ is an isomorphism. Then
$g\phi$ factors through $g'$ since $g'$ is a right
$\cal{D}$-approximation of $Y'$, and we may write $g\phi=hg'$ for
some $h: M\lra M'$. Similarly, there is a homomorphism $h': M'\lra
M$ such that $g'\phi^{-1}=h'g$. Thus
$hh'g=hg'\phi^{-1}=g\phi\phi^{-1}=g$ and
$h'hg'=h'g\phi=g'\phi^{-1}\phi=g'$. Since both $g$ and $g'$ are
right minimal, the morphisms $hh'$ and $h'h$ are isomorphisms. It
follows easily that $h$ itself is an isomorphism. Since $f'$ is a
kernel of $g'$ and since $f$ is a kernel of $g$, there is a morphism
$k: X\lra X'$ and a morphism $k': X'\lra X$ such that $kf'=fh$ and
$k'f=f'h^{-1}$. Thus $kk'f=kf'h^{-1}=fhh^{-1}=f$. It follows that
$kk'=1_{X}$ since $f$ is a monomorphism. Similarly, we have
$k'k=1_{X'}$. Hence $k$ is an isomorphism and the two sequences are
isomorphic. Similarly, the other statements in (2) can be proved.
$\square$

\medskip
To get an almost $\cal D$-split sequence, we may use the following
proposition. First, we introduce some notations. Let $\cal D$ be a
full subcategory of a category $\cal C$. An object $C$ in $\cal C$
is said to be \emph{generated} (respectively, \emph{co-generated})
by $\cal D$ if there is an epimorphism $D\lra C$ (respectively,
monomorphism $C\lra D$) with $D\in {\cal D}$. We denote by
${\mathscr F}({\cal D})$ the full subcategory of $\cal C$ consisting
of all objects $C\in \cal C$ generated by $\cal D$, and by
${\mathscr S}(\cal D)$ the full subcategory of $\cal C$ consisting
of all objects $C\in \cal C$ co-generated by $\cal D$.

\begin{Prop} Suppose $A$ is a ring with identity and $\cal C$ =
$A$-\emph{Mod}. Let $\cal D$ be a  full subcategory of $\cal C$. We
define ${\mathscr X}({\cal D})=\{ X\in {\cal C}\mid
\emph{Ext}_A^1(X, {\cal D})=0\}$ and ${\mathscr Y}({\cal D})=\{Y\in
{\cal C}\mid \emph{Ext}^1_A({\cal D},Y)=0\}$.

$(1)$ If $\cal D$ is contravariantly finite in $\cal C$, then, for
any $A$-module $Y\in {\mathscr F}({\cal D})\cap {\mathscr X}({\cal
D})$, there is an almost $\cal D$-split sequence $0\lra X\lra D\lra
Y\lra 0$ in $\cal C$.

$(2)$ If $\cal D$ is covariantly finite in $\cal C$, then, for any
$A$-module $X\in {\mathscr S}({\cal D})\cap {\mathscr Y}({\cal D})$,
there is an almost $\cal D$-split sequence $0\lra X\lra D\lra Y\lra
0$ in $\cal C$.
\end{Prop}

{\it Proof.} (1) Since $Y$ is generated by $\cal D$, there is a
surjective right $\cal D$-approximation of $Y$, say $g: M\lra Y$
with $M\in {\cal D}$. Let $X$ be the kernel of $g$. Then it follows
from the exact sequence $0\lra X\lra M\lra Y\lra 0$ that the
sequence Hom$_A(M,D')\lra \mbox{Hom}_A(X,D')\lra 0$ is exact since
$Y\in {\mathscr X}({\cal D})$. This implies that the homomorphism
$X\lra M$ is a left $\cal D$-approximation of $X$. Thus we get an
almost $\cal D$-split sequence in $\cal C$. (2) can be proved
analogously. $\square$

\medskip
Our main purpose of introducing almost $\cal D$-split sequences is
to construct derived equivalences between the endomorphism algebras
of objects appearing in almost $\cal D$-split sequences. The
following lemma is useful in our discussions.

\begin{Lem}
Let ${\cal C}$ be an additive category and $M$ an object in $\cal
C$. Suppose
$$X\stackrel{f}{\lra}M_n\lra\cdots\lra M_2\stackrel{t}{\lra} M_1\stackrel{g}{\lra}Y $$
is a (not necessarily exact) sequence of morphisms in $\cal C$ with
$M_i\in\add(M)$ satisfying the following conditions:

$(1)$ The morphism $f: X\lra M_n$ is a left $\add(M)$-approximation
of $X$, and the morphism $g: M_1\longrightarrow Y$ is a right
$\add(M)$-approximation of $Y$;

$(2)$ Put $V=M\oplus X$ and $W=M\oplus Y$. There are two induced
exact sequences
$$0\lra \Hom_{\cal C}(V, X)\stackrel{f_*}{\lra}\Hom_{\cal C}(V, M_n)\ra\cdots\ra
\Hom_{\cal C}(V,M_1)\stackrel{g_*}{\lra}\Hom_{\cal C}(V, Y), $$
$$0\lra \Hom_{\cal C}(Y, W)\stackrel{g^*}{\lra}\Hom_{\cal C}(M_1,W)\ra\cdots\ra
\Hom_{\cal C}(M_n,W)\stackrel{f^*}{\lra}\Hom_{\cal C}(X, W). $$ Then
the endomorphism rings $\End_{\cal C}(M\oplus X)$ and $\End_{\cal
C}(M\oplus Y)$ are derived-equivalent via a tilting module of
projective dimension at most $n$. \label{3.1}
\end{Lem}

{\it Proof.} Let $\Lambda$ be the endomorphism ring of $V$, and let
$T$ be the cokernel of the map $[t\,\,0]_*: \Hom_{\cal C}(V,
M_2)\longrightarrow \Hom_{\cal C}(V, M_1\oplus M)$. Then, by (2), we
have an exact sequence of $\Lambda$-modules:
$$ 0\ra \Hom_{\cal C}(V, X) \ra\Hom_{\cal C}(V, M_n)\ra \\
\cdots\ra \Hom_{\cal C}(V, M_2)\ra \Hom_{\cal C}(V, M_1\oplus M)\ra
T\ra 0. \quad (*)$$ Note that all the $\Lambda$-modules appearing in
the above exact sequence are finitely generated. Applying
$\Hom_{\Lambda}(-, \Hom_{\cal C}(V, M))$ to this sequence, we get a
sequence which is isomorphic to the following sequence
$$\begin{array}{rl}0\lra \Hom_{\Lambda}(T, \Hom_{\cal C}(V, M))& \lra\Hom_{\cal C}(M_1\oplus M, M)
\lra\Hom_{\cal C}(M_2, M)\lra \\
\\ & \cdots\lra \Hom_{\cal C}(M_n, M)\stackrel{f^*}{\lra} \Hom_{\cal
C}(X, M)\lra 0. \end{array}$$ By the second exact sequence in (2)
and the fact that $f$ is a left $\add(M)$-approximation of $X$, we
see that this sequence is exact. It follows that
$\Ext^i_{\Lambda}(T, \Hom_{\cal C}(V, M))=0$ for all $i>0$. Hence
$\Ext^i_{\Lambda}(T, \Hom_{\cal C}(V, M'))=0$ for all $i>0$ and
$M'\in\add(M)$. Thus, by applying $\Hom_{\Lambda}(T, -)$ to the
exact sequence $(*)$, we get
$\Ext^i_{\Lambda}(T,T)\simeq\Ext^{i+n}_{\Lambda}(T, \Hom_{\cal C}(V,
X))$ for all $i>0$. But $\Ext^{i+n}_{\Lambda}(T, \Hom_{\cal C}(V,
X))=0$ for all $i>0$ since the projective dimension of $T$ is at
most $n$. Thus Ext$^i_{\Lambda}(T,T)=0$ for all $i> 0$. Also, it
follows from the exact sequence $(*)$ that the following sequence
$$0\ra \Hom_{\cal C}(V, X\oplus M){\ra}\Hom_{\cal C}(V, M_n\oplus M)
\ra\cdots\ra \Hom_{\cal C}(V, M_2)\ra \Hom_{\cal C}(V, M_1\oplus
M)\ra T\ra 0$$ is exact, where $\Hom_{\cal C}(V, X\oplus M)$ is just
$\Lambda$ and the other terms are in $\add(T)$. Thus $T$ is a
tilting $\Lambda$-module of projective dimension at most $n$.

Next, we show that $\End_{\Lambda}(T)$ and $\End_{\cal C}(W)$ are
isomorphic. If $n=1$, we set $V'=X$ and $a=[f,0]: V'\lra M_1\oplus
M$. For $n\ge 2$, we set $V'=M_2$ and $a=[t,0]: V'\lra M_1\oplus M$.
Let $u: V'\lra V'$ and $v: M_1\oplus M\lra M_1\oplus M$ be two
morphisms in $\cal C$. The morphism pair $(u, v)$ is an endomorphism
of the sequence $V'\lra M_1\oplus M$ if $ua=av$. Let $E$ be the
endomorphism ring of the sequence $V'\lra M_1\oplus M$. Let $I$ be
the subset of $E$ consisting of those endomorphisms $(u, v)$ such
that there exists $h: M_1\oplus M\lra V'$ with $ha=v$. It is easy to
check that $I$ is an ideal of $E$. We shall show that $\End_{\cal
C}(W)$ is isomorphic to the quotient ring $E/I$. Let $b$ be the
morphism $\left[{{g}\atop 0}\,\,{{0}\atop{{\rm id}}}\right]:
M_1\oplus M\lra W$. Then, by the second exact sequence of the
condition (2), we have an exact sequence
$$0\lra\Hom_{\cal C}(W, W)\stackrel{b^*}{\lra}
\Hom_{\cal C}(M_1\oplus M, W)\stackrel{a^*}{\lra} \Hom_{\cal C}(V',
W). \hspace{1cm} (**)$$ By considering the image of ${\rm id}_W$
under the composition $b^*a^*$, we have $ab=0$. Thus, for each $(u,
v)\in E$, we have $avb=uab=0$, which means that $vb$ is in the
kernel of $a^*$. Therefore, there is a unique map $q: W\rightarrow
W$ such that $bq=vb$. Now, we define $\eta: E\rightarrow \End_{\cal
C}(W)$ by sending $(u, v)$ to $q$. Then $\eta$ is clearly a ring
homomorphism. We claim that $\eta$ is surjective. Indeed, since $g$
is a right $\add(M)$-approximation of $Y$, it is easy to check that
the map $b$ is a right $\add(M)$-approximation of $W$. Let $q$ be an
endomorphism of $W$. Then there is a morphism $v: M_1\oplus
M\longrightarrow M_1\oplus M$ such that $vb=bq$. By the first exact
sequence in $(2)$, we have the following exact sequence:
$$\Hom_{\cal C}(V', V')\stackrel{a_*}{\lra}\Hom_{\cal C}(V', M_1\oplus M)
\stackrel{b_*}{\lra}\Hom_{\cal C}(V', W).$$ It follows from
$avb=abq=0$ that $av$ is in the kernel of $b_*$ and there is a map
$u: V'\lra V'$ such that $ua=av$. This implies that $(u, v)$ is in
$E$ and $\eta(u,v)=q$. Hence $\eta$ is surjective.

Now, we determine the kernel of $\eta$. Note that, by the first
exact sequence in $(2)$, we have an exact sequence
$$\Hom_{\cal C}(M_1\oplus M, V')\stackrel{a_*}{\lra}\Hom_{\cal C}(M_1\oplus M, M_1\oplus M)
\stackrel{b_*}{\lra}\Hom_{\cal C}(M_1\oplus M, W). $$ Now, suppose
$(u, v)$ is in the kernel of $\eta$. Then $vb=0$, which means that
$v$ is in the kernel of $b_*$. Hence there is a map $h: M_1\oplus
M\lra V'$ such that $ha=v$. This implies $(u,v)\in I$. On the other
hand, if $(u, v)\in I$ and if $\eta$ sends $(u,v)$ to $q$, then
$bq=vb=hab=0$ and $q$ is in the kernel of $b^*$. By the exact
sequence $(**)$, we have $q=0$. Hence $I$ is the kernel of $\eta$,
and therefore $E/I\simeq \End_{\cal C}(W)$.

Let $\overline{E}$ be the endomorphism ring of  the following
complex of $\Lambda$-modules:
$$\Hom_{\cal C}(V, V')\stackrel{a_*}{\lra}\Hom_{\cal C}(V, M_1\oplus M),$$
and $\overline{I}$ the ideal of $\overline{E}$ consisting of those
$(\overline{u},\overline{v})$ such that
$\overline{h}a_*=\overline{v}$ for some $\overline{h}: \Hom_{\cal
C}(V, M_1\oplus M)\lra \Hom_{\cal C}(V, V')$. Similarly, we can show
that $\End_{\Lambda}(T)$ is isomorphic to
$\overline{E}/\overline{I}$. Finally, the natural map $e: E\lra
\overline{E}$, which sends $(u, v)$ to $(u_*, v_*)$, is clearly an
isomorphism of rings and induces an isomorphism from the ring $E/I$
to the ring $\overline{E}/\overline{I}$. Thus $\End_{\Lambda}(T)$
and $\End_{\cal C}(W)$ are isomorphic. The proof is completed.
$\square$

\medskip
{\it Remarks.} (1) For an Auslander-Reiten sequence $0\ra X\ra M\ra
Y\ra 0$ in $A$-mod with $A$ an Artin algebra, the proof that
End$(_AT)$ of the tilting module $T$ defined in Lemma \ref{3.1} is
isomorphic to $\End_A(M\oplus Y)$ can be carried out very easily.

(2) From the proof of Lemma \ref{3.1} we see that if we replace the
second exact sequence in (2) by the following two exact sequences
$$0\lra \Hom_{\cal C}(Y, M)\stackrel{g^*}{\lra}\Hom_{\cal C}(M_1,M)\ra\cdots\ra
\Hom_{\cal C}(M_n,M)\stackrel{f^*}{\lra}\Hom_{\cal C}(X, M), $$
$$0\lra \Hom_{\cal C}(Y,Y)\stackrel{g^*}{\lra} \Hom_{\cal C}(M_1,Y)\stackrel{t^*}{\lra} \Hom_{\cal C}(M_2,Y), $$
then Lemma \ref{3.1} still holds true. (Here $M_2=X$ if $n=1$.)
However, in most of cases that we are interested in, the second
exact sequence in (2) does exist.

(3) A special case of Lemma \ref{3.1} is the $n$-almost split
sequences in a maximal $(n-1)$-orthogonal subcategory studied in
\cite{iyama}. Let $A$ be a finite-dimensional algebra over a field.
Suppose $\cal C$ is a functorially finite and full subcategory of
$A$-mod. Recall that $\cal C$ is called a \emph{maximal}
$(n-1)$-\emph{orthogonal} subcategory if $\Ext^i_A(X,Y)=0$ for all
$X,Y\in {\cal C}$ and all $0<i\le n-1$, and ${\cal C}={\cal C}\cap
\{X\in A\mbox{-mod} \mid \Ext^i_A(C,X)=0 \mbox{\; for\;} C\in {\cal
C} \mbox{\; and\;} 0<i\le n-1\}$ = ${\cal C}\cap \{Y\in A\mbox{-mod}
\mid \Ext^i_A(Y,C)=0 \mbox{\; for\;} C\in {\cal C} \mbox{\; and\;}
0<i\le n-1\}$. In \cite{iyama}. It is shown that, for any
non-projective indecomposable $X$ in $\cal C$ (respectively,
non-injective indecomposable $Y$ in $\cal C$), there is an exact
sequence
$$ (*)\quad  0\ra Y\stackrel{f_n}{\lra}C_{n-1}\stackrel{f_{n-1}}{\lra}\cdots
\stackrel{f_1}{\lra}C_0\stackrel{f_0}{\lra} X\ra 0$$ with $C_j\in
{\cal C}$ and $f_j$ being radical maps such that the following
induced sequences are exact on $\cal C$:
$$  0\lra {\cal C}(-,Y)\lra {\cal C}(-,C_{n-1})\lra\cdots
\lra {\cal C}(-,C_0)\lra \mbox{rad}_{\cal C}(-,X)\lra 0,$$
$$  0\lra {\cal C}(X,-)\lra {\cal C}(C_0,-)\lra\cdots
\lra {\cal C}(C_{n-1},-)\lra \mbox{rad}_{\cal C}(Y,-)\lra 0,$$ where
rad$_{\cal C}$ stands for the Jacobson radical of the category $\cal
C$. Note also that $f_0$ is a minimal right almost split morphism
and that $f_n$ is a minimal left almost split morphism. The sequence
($*$) is called an $n$-almost split sequence in \cite{iyama}.

\medskip
With Lemma \ref{3.1} in mind, now we can show the significance of an
almost $\cal D$-split sequence for constructing derived equivalences
by the following result.

\begin{Theo}\label{exchangeTh}
Let $\cal C$ be an additive category and $M$ an object in $\cal C$.
Suppose
$$X\stackrel{f}{\lra}M'\stackrel{g}{\lra}Y $$
is an almost add$(M)$-split sequence in $\cal C$. Then the
endomorphism ring $\End_{\cal C}(M\oplus X)$ of $M\oplus X$ and the
endomorphism ring $\End_{\cal C}(M\oplus Y)$ of $M\oplus Y $ are
derived-equivalent.
\end{Theo}

\medskip
{\it Proof.} Let $V=M\oplus X$ and $W=M\oplus Y$. We shall verify
the conditions of Lemma \ref{3.1} for $n=1$. By the definition of an
almost $\cal D$-split sequence, we see immediately that the
condition (1) in Lemma \ref{3.1} is satisfied, while the condition
(2) in Lemma \ref{3.1} is implied by the condition (3) in Definition
\ref{def1}: In fact, by applying $\Hom_{\cal C}(V,-)$ to the above
sequence, we get a complex of abelian groups
$$(*) \quad 0\lra \Hom_{\cal C}(V,X)\stackrel{(-,f)}{\lra}\Hom_{\cal
C}(V,M')\stackrel{(-,g)}{\lra} \Hom_{\cal C}(V,Y).$$ Since $f$ is a
monomorphism, the map $(-,f)$ is injective. Clearly, the image of
the map $(-,f)$ is contained in the kernel of the map $(-,g)$. Since
$f$ is a kernel of $g$, it is easy to see that the kernel of $(-,g)$
is equal to the image of $(-,f)$. Thus ($*$) is exact. Similarly, we
see that the sequence
$$0\lra \Hom_{\cal
C}(Y,W)\stackrel{(g,-)}{\lra}\Hom_{\cal
C}(M',W)\stackrel{(f,-)}{\lra} \Hom_{\cal C}(X,W)$$ is exact. Thus
Theorem \ref{exchangeTh} follows from Lemma \ref{3.1} if we take
$n=1$. $\square$

\medskip
In Theorem \ref{exchangeTh}, the two rings $\End_{\cal C}(M\oplus
X)$ and $\End_{\cal C}(M\oplus Y)$ are linked by a tilting module of
projective dimension at most $1$. This is precisely the case of
classic tilting module. Thus there is a nice linkage between the
torsion theory defined by the tilting module in $\End_{\cal
C}(M\oplus X)$-mod and the one in $\End_{\cal C}(M\oplus Y)$-mod.
For more details we refer to \cite{BB} and \cite{hr}.

In the following, we deduce some consequences of Theorem
\ref{exchangeTh}. Since an Auslander-Reiten sequence can be viewed
as an almost $\cal D$-split sequence, as explained in Example (b),
we have the following corollary.

\begin{Koro}
Let $A$ be an Artin algebra, and let $0\ra X\ra M\ra Y\ra 0$ be an
Auslander-Reiten sequence in $A\modcat$. Suppose $N$ is an
$A$-module in $A$-mod such that neither $X$ nor $Y$ belongs to
$\add(N)$. Then $\End_A(N\oplus M\oplus X)$ is derived-equivalent to
$\End_A(N\oplus M\oplus Y)$. In particular, $\End_A(M\oplus X)$ and
$\End_A(M\oplus Y)$ are derived-equivalent.\label{ars}
\end{Koro}

As another consequence of Theorem \ref{exchangeTh}, we have the
following corollary.

\begin{Koro}\label{selinjsysygy} Let $A$ be an Artin algebra and $X$ a torsion-less $A$-module, that is,
$X$ is a submodule of a projective module in $A\emph{-mod}$. If $f:
X\ra P$ is a left $\add(_AA)$-approximation of $X$, then
$\End_A(A\oplus X)$ and $\End(_AA\oplus \emph{coker}(f))$ are
derived-equivalent. In particular, if $A$ is a self-injective Artin
algebra, then, for any $X$ in $A\emph{-mod}$, the algebras
$\End_A(A\oplus X)$ and $\End_A(A\oplus \Omega (X))$ are
derived-equivalent via a tilting module. \label{cor3.6}
\end{Koro}

{\it Proof.} Note that $f$ is injective. Thus the short exact
sequence $$0\lra X\stackrel{f}{\lra} P\lra \mbox{coker}(f)\lra 0$$
is an almost add$(_AA)$-split sequence in $A$-mod. By Theorem
\ref{exchangeTh}, the corollary follows. $\square$

\medskip As a consequence of Corollary \ref{selinjsysygy}, we
get the following corollary.

\begin{Koro}
Let $A$ be a self-injective Artin algebra and $X$ an $A$-module.
Then the algebras $\End_A(A\oplus X)$ and $\End_A(A\oplus\tau X)$
are derived-equivalent, where $\tau$ stands for the Auslander-Reiten
translation. Thus, for all $n$, the algebras $\End_A(A\oplus
\tau^nX)$ are derived-equivalent. \label{cor3.7}
\end{Koro}

{\it Proof.} Let $\nu$ be the Nakayama functor $D\Hom_A(-,A)$. It is
known that if $A$ is self-injective then $\tau \simeq \nu\Omega^2$,
$\nu(A)=A$ and the Nakayama functor is an equivalence from $A$-mod
to itself. Since the algebra $\End_A(A\oplus \tau X)$ is isomorphic
to the algebra $\End_A(A\oplus\Omega^2(X))$, the corollary follows
immediately from Corollary \ref{selinjsysygy}. $\square$

\medskip
{\it Remark.} If $A$ is a finite-dimensional self-injective algebra,
then, for any $A$-module $X$, it was shown in \cite[Corollary
1.2]{LiuXi3} that the algebras $\End_A(A\oplus X)$, $\End_A(A\oplus
\Omega (X))$ and $\End_A(A\oplus \tau X)$ are stably equivalent of
Morita type. Thus they are both derived-equivalent and stably
equivalent of Morita type. For further information on stably
equivalences of Morita type for general finite-dimensional algebras,
we refer the reader to \cite{LiuXi1, LiuXi2, LiuXi3, XiAdjointtype}
and the references therein.

Now, we point out the following consequence of Theorem
\ref{exchangeTh}: if $0\ra X\ra M'\ra Y\ra 0$ is an almost $\cal
D$-split sequence in $A$-mod with $\cal D$ = add$(M)$ for an
$A$-module $M$, then $X$ and $Y$ have the same number of
non-isomorphic indecomposable direct summands which are not in
add$(M)$. This follows from the fact that a derived equivalence
preserves the number of non-isomorphic simple modules.

Many other invariants of derived equivalences can be used to study
the algebras End$_A(M\oplus X)$ and End$_A(M\oplus Y)$; for example,
End$_A(M\oplus X)$ has finite global dimension if and only if
End$_A(M\oplus Y)$ has finite global dimension. This follows from
the fact that derived equivalence preserves the finiteness of global
dimension. In fact, we have the following explicit formula by
tilting theory (see \cite{hr} and \cite[Proposition 3.4,
p.116]{HappelTri}, for example):

If $0\ra X\ra M'\ra Y\ra 0$ is an almost $\cal D$-split sequence in
$A$-mod with $\cal D$ = add$(M)$ for an $A$-module $M$ in $A$-mod,
then
$$\mbox{gl.dim}(\End_{\cal C}(M\oplus X))-1\le \mbox{gl.dim}(\End_{\cal
C}(M\oplus Y))\le \mbox{gl.dim}(\End_{\cal C}(M\oplus X))+1,$$ where
gl.dim$(A)$ stands for the global dimension of $A$. Note that the
global dimension of $\End_{\cal C}(M\oplus X))$ may be infinite (see
Example $3$ in Section \ref{sect7}). Concerning global dimensions
and Auslander-Reiten sequences, there is a related result which can
be found in \cite{XH}.

Note that if a derived equivalence between two rings $A$ and $B$ is
obtained from a tilting module $_AT$, that is, there exists a
tilting $A$-module $_AT$ such that $B\simeq \End_A(T)$,  then the
finitistic dimension of $A$ is finite if and only if the finitistic
dimension of $B$ is finite (see \cite{ha})\footnote{Recently, it is
shown that the finiteness of finitistic dimension is invariant under
an arbitrary derived equivalence.}. Recall that the \emph{finitistic
dimension} of an Artin algebra $A$, denoted by fin.dim$(A)$, is
defined to be the supremum of the projective dimensions of finitely
generated $A$-modules of finite projective dimension. The finitistic
dimension conjecture states that fin.dim$(A)$ should be finite for
any Artin algebra $A$. This conjecture has closely been related to
many other homological conjectures in the representation theory of
algebras. For some advances and further information on the
finitistic dimension conjecture, we may refer the reader to the
recent paper \cite{xx} and the references therein.

Thus we have the following corollary.

\begin{Koro} Let $\cal C$ be an additive category and $M$ an object in $\cal C$. Suppose
$$X\stackrel{f}{\longrightarrow}M'\stackrel{g}{\longrightarrow}Y
$$ is an almost \emph{add}$(M)$-split sequence in
$\cal C$. Then the finitistic dimension of $\End_{\cal C}(M\oplus
X)$ is finite if and only if the finitistic dimension of $\End_{\cal
C}(M\oplus Y)$ is finite.
\end{Koro}

\medskip
If $A$ is an Artin $R$-algebra over a commutative Artin ring
$R$ and $M$ is an $A$-bimodule, then $A\ltimes M$, the trivial
extension of $A$ by $M$ is the $R$-algebra whose underlying
$R$-module is $A\oplus M$, with multiplication given by

$$ (\lambda,m)(\lambda', m')= (\lambda\lambda',\lambda m'+
m\lambda')$$ for $\lambda,\lambda'\in A$, and $m, m'\in M$. It is
shown in  \cite{RickDstable} that if $A$ and $B$ are
finite-dimensional algebras over a field $k$ that are
derived-equivalent, then $A\ltimes D(A)$ is derived-equivalent to
$B\ltimes D(B)$, where $D=\Hom_k(-,k)$. Note that $A\ltimes D(A)$ is
a self-injective algebra and that a derived equivalence between two
self-injective algebras implies a stable equivalence of Morita type
between them by \cite{RickDstable}. It is  known in \cite{Xi1} that
a stable equivalence of Morita type preserves representation
dimension (see \cite{Auslander} for definition). Hence we have the
following corollary.

\begin{Koro} Let $\Lambda$ be a finite-dimensional algebra over a field $k$ and
$M$ a $\Lambda$-module in $\Lambda\modcat$. Suppose
$$X\stackrel{f}{\longrightarrow}M'\stackrel{g}{\longrightarrow}Y
$$ is an almost \emph{add}$(M)$-split sequence in
$\Lambda\modcat$, and let $A=\End_{\Lambda}(X\oplus M)$ and
$B=\End_{\Lambda}(M\oplus Y)$. Then $A\ltimes D(A)$ is
derived-equivalent to $B\ltimes D(B)$. In particular, the
representation dimensions of $A\ltimes D(A)$ and $B\ltimes D(B)$ are
equal.
\end{Koro}

\medskip
In the following, we consider several generalizations of Corollary
\ref{ars}, namely we deal with the case of a finite family of
Auslander-Reiten sequences.

\begin{Koro}\label{corAR-chain}
Let $A$ be an Artin algebra, and let $0\lra X_i\lra M_i\lra
X_{i-1}\lra 0$ be an Auslander-Reiten sequence in $\modcat{A}$ for
$i=1,2,\cdots, n$. Let $M=\bigoplus_{i=1}^nM_i$.  Then
$\End_A(M\oplus X_n)$ and $\End_A(M\oplus X_0)$ are
derived-equivalent via a tilting module $T$ of projective dimension
at most $n$.
\end{Koro}

{\it Proof.} First, we suppose $X_n\in\add(M)$. Then there is an
$M_i$ such that $X_n$ is a direct summand of $M_i$, and therefore
there is an irreducible map from $X_i$ to $X_n$. It follows that
there is an irreducible map from $X_0=\tau^{-i}X_i$ to
$X_{n-i}=\tau^{-i}X_n$, where $\tau$ stands for the Auslander-Reiten
translation. Thus $X_0$ is a direct summand of $M_{n-i+1}$, which
implies $X_0\in\add(M)$. Hence $\add(M\oplus
X_n)=\add(M)=\add(M\oplus X_0)$. Consequently, the algebras
$\End_A(M\oplus X_n)$ and $\End_A(M\oplus X_0)$ are Morita
equivalent. Thus $\End_A(M\oplus X_n)$ and $\End_A(M\oplus X_0)$
are, of course, derived-equivalent via a (projective) tilting
module.

Next, we assume $X_n\not\in\add(M)$. In this case, we claim that
there is no integer $i\in\{0,1,\cdots,n\}$ such that
$X_i\in\add(M)$. If $X_0\in\add(M)$, then there is an $M_i, 1\leq
i\leq n,$ such that $X_0$ is a direct summand of $M_i$. Thus there
is an irreducible map from $X_i$ to $X_0$. By applying the
Auslander-Reiten translation, we see that there is an irreducible
map from $X_n=\tau^{n-i}X_i$ to $X_{n-i}=\tau^{n-i}X_0$. Hence $X_n$
is a direct summand of $M_{n-i+1}$, that is, $X_n$ is in $\add(M)$.
This is a contradiction and shows that $X_0$ does not belong to
$\add(M)$. Suppose $X_i\in\add(M)$ for some $0<i<n$. Then there is
an integer $j \in \{1,2, \cdots, n\}$ such that $X_i$ is a direct
summand of $M_j$. Clearly, $i\neq j$, and there is an irreducible
map from $X_i$ to $X_{j-1}$. On the one hand, if $i>j$, then there
is an irreducible map from $X_n=\tau^{n-i}X_i$ to
$X_{n-i+j-1}=\tau^{n-i}X_{j-1}$. This implies that $X_n$ is a direct
summand of $M_{n-i+j}$, which is a contradiction. On the other hand,
if $i<j$, then there is an irreducible map from $X_0=\tau^{-i}X_i$
to $X_{j-1-i}=\tau^{-i}X_{j-1}$. It follows that $X_0$ is a direct
summand of $M_{j-i}$. This is again a contradiction. Hence there is
no $X_i$ belonging to $\add(M)$.

Now let $m$ be the minimal integer in $\{0,1,\cdots, n\}$ such that
$X_n\simeq X_m$. If $m=0$, then $\add(M\oplus X_n)=\add(M\oplus
X_0)$. This means that the endomorphism algebras $\End_A(M\oplus
X_n)$ and $\End_A(M\oplus X_0)$ are Morita equivalent. Now we assume
$m>0$. Then the $A$-modules $X_0,X_1,\cdots,X_m$ are pairwise
non-isomorphic. We consider the sequence
$$X_m\longrightarrow M_m\longrightarrow\cdots\longrightarrow M_1\longrightarrow X_0.$$
Since $X_m\not\in\add(M)$, every homomorphism from $X_m$ to $M$
factors through the map $X_m\lra M_m$ in the Auslander-Reiten
sequence starting at $X_m$. This means that the map
$X_m\longrightarrow M_m$ is a left $\add(M)$-approximation of $X_m$.
Similarly, the map $M_1\lra X_0$ is a right $\add(M)$-approximation
of $X_0$. Let $V=M\oplus X_m$. Then  $X_i\not\in\add(V)$ for all
$i=0,1,\cdots, m-1$. It follows that we have exact sequences
$$0\lra \Hom_A(V, X_i)\lra \Hom_A(V, M_i)\lra \Hom_A(V, X_{i-1})\lra 0$$
for $i=1,\cdots, m$. Connecting the above exact sequences, we get an
exact sequence
$$0\lra\Hom_A(V, X_m)\lra\Hom_A(V, M_m)\lra\cdots\lra\Hom_A(V, M_1)\lra \Hom_A(V, X_0).$$
This gives the first exact sequence in Lemma \ref{3.1}(2). The
second exact sequence in Lemma \ref{3.1}(2) can be obtained
similarly. Thus Corollary \ref{corAR-chain} follows immediately from
Lemma \ref{3.1}. $\square$

\medskip
{\it Remark.} In Corollary \ref{corAR-chain}, if $X_n\not\in
\add(M)$ and $X_0, X_1,\cdots, X_n$ are pairwise non-isomorphic,
then the tilting $\End(X\oplus M)$-module $T$ defined in Lemma
\ref{3.1} has projective dimension $n$. Note that we always have
gl.dim$(\End_A(X\oplus M))-n\le$ gl.dim$(\End_A(M\oplus Y))\le $
gl.dim$(\End_A(X\oplus M))+n.$

The following is another type of generalization of Corollary
\ref{ars}.

\begin{Prop} Let $A$ be an Artin algebra.

$(1)$ Suppose $0\lra X_i\lra M_i\lra Y_i\lra 0 $ is an
Auslander-Reiten sequence for $i=1, 2, \cdots, n$.  Let
$X=\bigoplus_i X_i, M=\bigoplus_iM_i$ and $Y=\bigoplus_iY_i$. If
$\add(X)\cap \add(M) = 0$ = $\add(M)\cap \add(Y)$, then
$\End_A(X\oplus M)$ and $\End_A(M\oplus Y)$ are derived-equivalent.

$(2)$ Suppose $0\lra X_1\lra X_2\oplus M_1\lra Y_1\lra 0$ and $0\lra
X_2\lra Y_1\oplus M_2\lra Y_2\lra 0$ are two Auslander-Reiten
sequences such that neither $X_2$ is in $\add(M_1)$ nor $Y_1$ is in
$\add(M_2)$. If $X_1\not\in \add(Y_1\oplus M_2)($or equivalently,
$Y_2\not\in \add(X_2\oplus M_1))$, then $\End_A(X_1\oplus M_1\oplus
M_2)$ and $\End_A( M_1\oplus M_2\oplus Y_2)$ are derived-equivalent.
\label{generalization}
\end{Prop}

{\it Proof.} (1) Under our assumption, the exact sequence $ 0\lra
X\lra M\lra Y\lra 0$ is an almost add$(M)$-split sequence in
$A$-mod. Therefore (1) follows from Theorem \ref{exchangeTh}.

(2) There is an exact sequence

$$(*) \quad 0\lra X_1\lra M_1\oplus M_2\lra Y_2\lra 0,$$
which can be constructed by the given two  Auslander-Reiten
sequences. Clearly, $X_1\not\in \add(X_2\oplus M_1)$ since
Auslander-Reiten quiver has no loops. By assumption, we see
$X_1\not\in \add(M_1\oplus M_2)$. Hence we can verify that the
morphism $X_1\lra M_1\oplus M_2$ in ($*$) is a left $\add(M_1\oplus
M_2)$-approximation of $X_1$. Similarly, we can see that the
morphism $ M_1\oplus M_2\lra Y_2$ in ($*$) is a right
$\add(M_1\oplus M_2)$-approximation of $Y_2$. Thus ($*$) is an
almost $\add(M_1\oplus M_2)$-split sequence in $A$-mod, and
therefore the conclusion (2) follows from Theorem \ref{exchangeTh}.
$\square$

\medskip
{\it Remark.}  Usually, given two Auslander-Reiten sequences $0\ra
X_i\ra M_i\ra Y_i\ra 0$ ($1\le i\le 2$), we cannot get a derived
equivalence between $\End_A(X_1\oplus X_2\oplus M_1\oplus M_2)$ and
$\End_A( M_1\oplus M_2\oplus Y_1\oplus Y_2)$. For a counterexample,
we refer the reader to Example 3 in the last section.

Now, we mention that, for an $n$-almost split sequence studied in
\cite{iyama}, we have a statement similar to Corollary
\ref{corAR-chain}.

\begin{Prop} Let $\cal C$ be a maximal $(n-1)$-orthogonal subcategory
of $A$-\emph{mod} with $A$ a finite-dimensional algebra over a field
$(n\ge 1)$. Suppose $X$ and $Y$ are two indecomposable $A$-modules
in $\cal C$ such that the sequence
$$0\lra X\stackrel{f}{\lra} M_n\stackrel{t_n}{\lra}M_{n-1}\lra\cdots
\lra M_2\stackrel{t_2}{\lra} M_1\stackrel{g}{\lra} Y\lra 0$$ is an
$n$-almost split sequence in $\cal C$. Then $\End_A(X\oplus
\bigoplus_{i=1}^nM_i)$ and $\End_A(\bigoplus_{i=1}^nM_i\oplus Y)$
are derived-equivalent. \label{n-ars}
\end{Prop}

{\it Proof.} Let $M:=\bigoplus_{i=1}^{n}M_i$. Suppose $Y$ is a
direct summand of some $M_i$. Then there is a canonical projection
$\pi: M_i\lra Y$. Let $t_1=g$ and $t_{n+1}=f$. We observe that all
homomorphisms $t_1,\cdots,t_{n+1}$ are radical maps by the
definition of an $n$-almost split sequence. Hence the composition
$t_{i+1}\pi$ can not be a split epimorphism and consequently factors
through $t_1=g$, that is, $t_{i+1}\pi=u_1g$ for a homomorphism
$u_1:M_{i+1}\lra M_1$. First, we assume that $i\neq n$. Then
$t_{i+2}u_1g=t_{i+2}t_{i+1}\pi=0$. By \cite[Theorem 2.5.3]{iyama},
we have $t_{i+2}u_1=u_2t_2$ for a homomorphism $u_2: M_{i+2}\lra
M_2$. Similarly, we get a homomorphism $u_k: M_{i+k}\lra M_k$ such
that $t_{i+k}u_{k-1}=u_kt_k$ for $k=2,3,\cdots,n-i$. This allows us
to form the following commutative diagram:
$$\xymatrix{
 X\ar[r]^{f}\ar[d]^{u_{n-i+1}} & M_n\ar[d]^{u_{n-i}}\ar[r]^{t_n}
 & M_{n-1}\ar[d]^{u_{n-i-1}}\ar[r] & \cdots \ar[r] &
 M_{i+1}\ar[d]^{u_1}\ar[r]^{t_{i+1}} & M_i\ar[d]^{\pi}\ar[r]^{t_i} & M_{i-1}\\
M_{n-i+1}\ar[r]^{t_{n-i+1}} & M_{n-i}\ar[r]^{t_{n-i}} &
M_{n-i-1}\ar[r] & \cdots \ar[r] &
 M_{1}\ar[r]^{g} & Y.\\
}$$ Note that if $i=n$ then the above diagram still makes sense. We
claim that $u_{n-i+1}$ is a split monomorphism.  If this is not the
case, then the map $u_{n-i+1}$ factors through $f$. This means that
there is some map $h_n: M_n\lra M_{n-i+1}$ such that
$fh_n=u_{n-i+1}$. Then we have
$f(u_{n-i}-h_nt_{n-i+1})=fu_{n-i}-u_{n-i+1}t_{n-i+1}=0$. By
\cite[Theorem 2.5.3]{iyama}, there is some homomorphism $h_{n-1}:
M_{n-1}\lra M_{n-i}$ such that $t_nh_{n-1}=u_{n-i}-h_nt_{n-i+1}$,
that is, $u_{n-i}=t_nh_{n-1}+h_nt_{n-i+1}$. Similarly, we get $h_k:
M_k\lra M_{k-i+1}$ such that $u_{k-i+1}=h_{k+1}t_{k-i+2}+t_{k-i+1}h$
for $k=n-2,n-3,\cdots, i$. Thus
$t_{i+1}(\pi-h_ig)=t_{i+1}\pi-(u_{i+1}-h_{i+1}t_2)g=t_{i+1}\pi-u_{i+1}g=0$.
Hence $\pi-h_ig$ factors through $t_i$, say $\pi-h_ig=t_ih_{i-1}$.
Then $\pi=h_ig+t_ih_{i-1}$, which is a radical map since both $g$
and $t_i$ are radical maps. This is a contradiction. Hence $X$ is a
direct summand of $M_{n-i+1}$ and $\add(M\oplus
X)=\add(M)=\add(M\oplus Y)$. Thus, $\End_A(M\oplus X)$ and
$\End_A(M\oplus Y)$ are Morita equivalent.

Similarly, if $X$ is a direct summand of some $M_i$, then $Y$ is a
direct summand of $M_{n-i+1}$. It follows that $\End_A(M\oplus X)$
and $\End_A(M\oplus Y)$ are Morita equivalent.

Now we assume that neither $X$ nor $Y$ is a direct summand of $M$.
We use Lemma \ref{3.1} to show Proposition \ref{n-ars}. By a
property of an $n$-almost split sequence (see \cite[Theorem
2.5.3]{iyama}) and the fact that $X$ and $Y$ do not lie in add$(M)$,
we see that $f$ is a left add$(M)$-approximation of $X$ and $g$ is a
right add$(M)$-approximation of $Y$. It remains to check the
condition (2) in Lemma \ref{3.1}. However, it follows from
\cite[Theorem 2.5.3]{iyama} (see Remark (3) at the end of the proof
of Lemma \ref{3.1}) that we have two exact sequences
$$0\lra \Hom_A(V,X)\stackrel{(-,f)}{\lra}\Hom_A(V,M_n)\lra\cdots
\lra \Hom_A(V,M_1)\stackrel{(-,g)}{\lra} \Hom_A(V,Y),$$
$$0\lra \Hom_A(Y,W)\stackrel{(g,-)}{\lra} \Hom_A(M_1,W)\lra\cdots
\lra \Hom_A(M_n,W)\stackrel{(f,-)}{\lra} \Hom_A(X,W)$$ for
$V:=X\oplus M$ and $W:=M\oplus Y$. Thus the condition (2) in Lemma
\ref{3.1} is satisfied. Consequently, Proposition \ref{n-ars}
follows from Lemma \ref{3.1}. $\square$

\section{Auslander-Reiten sequences and BB-tilting modules\label{sect3+}}

In this section, we point out that, when we restrict our
consideration to Auslander-Reiten sequences, the tilting module
defining the derived equivalence in Theorem \ref{exchangeTh} is of
special form, namely it is a BB-tilting-module in the sense of
Brenner and Butler \cite{BB}. This shows that the tilting theory and
the Auslander-Reiten theory are so beautifully integrated with each
other. We first recall the BB-tilting-module procedure in \cite{BB},
and then give a generalization of a BB-tilting module, namely the
notion of an $n$-BB-tilting module.

Let $A$ be an Artin algebra and $S$ a non-injective simple
$A$-module with the following two properties: (a)
proj.dim$_A(\tau^{-1}S)$ $ \le 1$, and (b) Ext$^1_A(S,S)=0$. Here
$\tau^{-1}$ stands for the Auslander-Reiten translation TrD, and
proj.dim$_A(S)$ means the projective dimension of $S$. We denote the
projective cover of $S$ by $P(S)$, and assume that $_AA= P(S)\oplus
P$ such that there is not any direct summand of $P$ isomorphic to
$P(S)$. Let $T= \tau^{-1}S\oplus P$. It is well-known that $T$ is a
tilting module. Such a tilting module is called a BB-\emph{tilting
module}. In particular, if $S$ is a projective non-injective simple
module, then $T$ is automatically a BB-tilting module, this special
case was first studied in \cite{APR}, and the tilting module of this
form is called an APR-\emph{tilting module} in literature. Note that
if $S$ is a non-injective, projective simple $A$-module, then there
is an Auslander-Reiten sequence
$$  0\lra S\lra P'\lra \tau^{-1}S\lra 0$$
in $A$-mod with $P'$ projective.

\begin{Prop} Let $A$ be an Artin algebra, and let $0\lra X\stackrel{f}{\lra}
M\stackrel{g}{\lra} Y\lra 0$ be an Auslander-Reiten sequence in
$A\modcat$. We define $V:= M\oplus X$, $\Lambda =\End_A(V)$,
$W=M\oplus Y$ and $\Gamma=\End_A(W)$. Then the derived equivalence
between $\Lambda$ and $\Gamma$ in Theorem \ref{exchangeTh} is given
by a BB-tilting module. In particular, if the Auslander-Reiten
sequence
$$ 0\lra S\lra P'\lra \tau^{-1}S\lra 0 $$
defines an APR-tilting module $T:=P\oplus \tau^{-1}S$, then the
sequence is an almost $\add(P)$-split sequence in $A\modcat$ and the
derived equivalence between $A$ and $\End_A(T)$ in Theorem
\ref{exchangeTh} is given precisely by the APR-tilting module
$T:=P\oplus \tau^{-1}S$. \label{apr}
\end{Prop}

{\it Proof.}  From the Auslander-Reiten sequence we have the
following exact sequence
$$  0\ra \Hom_A(V,X)\ra\Hom_A(V,M)\stackrel{(-,g)}{\lra}\Hom_A(V,Y).$$
Let $L$ be the image of the map $(-,g)$. Then we have an exact
sequence
$$ (*)\quad 0\ra \Hom_A(V,X)\ra\Hom_A(V,M)\stackrel{(-,g)}{\lra}L\ra
0.$$ (This is a minimal projective presentation of the
$\Lambda$-module $L$). Let $T:= L\oplus \Hom_A(V,M)$. Then $T$ is
the tilting module which defines the derived equivalence in Theorem
\ref{exchangeTh}. We shall show that $T$ is a BB-tilting
$\Lambda$-module. To prove this, it is sufficient to show that $L$
is of the form $\tau^{-1}S$ for a simple $\Lambda$-module $S$.

If we apply $\Hom_{\Lambda}(-,\Lambda)$ to ($*$), then we get an
exact sequence of right $\Lambda$-modules:   $$
\Hom_{\Lambda}(\Hom_A(V,M),\Lambda)\lra
\Hom_{\Lambda}(\Hom_A(V,X),\Lambda)\lra \mbox{Tr}_{\Lambda}(L)\lra
0,$$ which is isomorphic to the following exact sequence
$$ \Hom_A(M,V)\stackrel{(f,-)}{\lra} \Hom_A(X,V)\lra \mbox{Tr}_{\Lambda}(L)\lra
0,$$ where Tr$_{\Lambda}$ stands for the transpose over $\Lambda$.
Note that the image of the map $(f,-)$ is the radical of the
indecomposable projective right $\Lambda$-module $\Hom_A(X,V)$. Thus
Tr$_{\Lambda}(L)$ is a simple right $\Lambda$-module, and
consequently, $\tau_{\Lambda}L$ is isomorphic to the socle $S$ of
the indecomposable injective $\Lambda$-module $D\Hom_A(X,V)$. Hence
$L\simeq \tau_{\Lambda}^{-1}S$. Since $X$ is not a direct summand of
$M$, we see that $\Ext^1_{\Lambda}(S,S)=0$. Thus $T$ is a BB-tilting
module. Note that if $X\not\simeq Y$ then $L\simeq \Hom_A(V,Y)$. In
case of an APR-tilting module, we can see that the given
Auslander-Reiten sequence is an almost $\add(P)$-split sequence.
Thus Proposition \ref{apr} follows. $\square$

\medskip
Now, we introduce the notion of an $n$-BB-tilting module: Let $A$ be
an Artin algebra. Recall that we denote by $\Omega^n$ the $n$-th
syzygy operator, and by $\Omega^{-n}$ the $n$-th co-syzygy operator.
As usual, $D$ is the duality of an Artin algebra. Suppose $S$ is a
simple $A$-module and $n$ is a positive integer. If $S$ satisfies
(a) $\Ext^j_A(D(A),S)=0$ for all $0\le j\le n-1$, and (b)
$\Ext^i_A(S,S)=0$ for all $1\le i\le n$, we say that $S$ defines an
$n$-BB-tilting module, and that the module $T :=
\tau^{-1}\Omega^{-n+1}(S)\oplus P$ is an $n$-BB-\emph{tilting
module}, where $P$ is the direct sum of all non-isomorphic
indecomposable projective $A$-modules which are not isomorphic to
$P(S)$, the projective cover of $S$.  Note that (a) implies that the
injective dimension of $S$ is at least $n$ and that the case $n=1$
is just the usual BB-tilting module. The terminology is adjudged by
the following lemma.

\begin{Lem} If $S$ defines an $n$-BB-tilting $A$-module, then $T:=\tau^{-1}\Omega^{-n+1}S\oplus P$
is a tilting module of projective dimension at most $n$.
\end{Lem}

{\it Proof.} Let $\nu$ be the Nakayama functor $D\Hom_A(-,{}_AA)$.
Suppose the sequence
$$0\lra S\lra \nu P_0\lra \nu P_1\lra \cdots\lra \nu P_n \lra \cdots
$$ is a minimal injective resolution of $S$ with all $P_i$
projective. Since $\Ext^i_A(D(A),S)=0$ for $0\le i\le n-1$, we have
the following exact sequence by applying $\Hom_A(D(A),-)$ to the
injective resolution:

$$0 \lra\Hom_A(D(A),S)\lra\Hom_A(D(A),\nu P_0)\lra\cdots \lra
\Hom_A(D(A),\nu P_n)\lra L\lra 0,$$ which is isomorphic to the
following exact sequence
$$0\lra 0\lra P_0\lra\cdots \lra P_n\lra L\lra 0. $$
This shows that $L\simeq \mbox{TrD}\Omega_A^{-n+1}(S)$ and the
projective dimension of $L$ is at most $n$. Moreover, we have the
following sequence:
$$ (*)\quad  0\lra \Hom_A(L,P)\lra \Hom_A(P_n,P)\lra \cdots \lra \Hom_A(P_0,P)\lra 0. $$
Since $\Hom_A(\nu P_j,\nu P)\simeq \Hom_A(P_j,P)$, we see that ($*$)
is isomorphic to the sequence
$$ 0\lra \Hom_A(L,P)\lra \Hom_A(\nu P_n,\nu P)\lra \cdots \lra \Hom_A(\nu P_0,\nu P)\lra 0, $$
which is exact because $\Hom_A(-,\nu P)$ is an exact functor. Note
that $\Hom_A(S,\nu P)=0$ by the definition of $P$. This shows that
$\Ext^i_A(L,P)=0$ for all $i>0$. Since $\Ext^i_A(S,S)=0$ for all
$1\le i\le n$, this means that $\nu P_0$ is not a direct summand of
$\nu P_i$ for $1\le i\le n$. Thus $P(S)$ is not a direct summand of
$P_i$ for $1\le i\le n$, that is, $P_i\in \add(P)$ for all $1\le
i\le n$. Now, if we apply $\Hom_A(L,-)$ to the projective resolution
of $L$, we get $\Ext^{n+i}_A(L,P_0)\simeq \Ext^i_A(L,L)$ for all
$i\ge 1$. Hence $\Ext^i_A(L,L)=0$ for all $i\ge 1$.

We note that $P_0 = P(S)$ and there is an exact sequence
$$ 0\lra A\lra P\oplus P_1\lra \cdots\lra L\lra 0.$$
Altogether, we have shown that $T$ is a tilting module of projective
dimension at most $n$. $\square$

\begin{Prop} $(1)$ Suppose $0\ra X_i\ra M_i\ra X_{i-1}\ra 0$ is an
Auslander-Reiten sequence in $A\modcat$ for $i=1,2,\cdots,n$. Let
$M=\bigoplus_{i=1}^n M_i$ and $V=M\oplus X_n$. If $X_n\not\in
\add(M)$ and if $X_0, X_1,\cdots, X_n$ are pairwise non-isomorphic,
then the tilting $\End_A(V)$-module $T:=\Hom_A(V,X_0)\oplus
\Hom_A(V,M)$ is an $n$-BB-tilting module.

$(2)$ Let $\cal C$ be a maximal $(n-1)$-orthogonal subcategory of
$A$-\emph{mod} with $A$ a finite-dimensional algebra over a field
$(n\ge 1)$. Suppose $X$ and $Y$ are two indecomposable $A$-modules
in $\cal C$ such that the sequence
$$0\lra X\stackrel{f}{\lra} M_n\stackrel{t_n}{\lra}M_{n-1}\lra\cdots
\lra M_2\stackrel{t_2}{\lra} M_1\stackrel{g}{\lra} Y\lra 0$$ is an
$n$-almost split sequence in $\cal C$. We define
$V=\bigoplus_{i=1}^n M_i\oplus X$, and $L$ to be the image of the
map $\Hom_A(V,g)$. If $X\not\in \add(\bigoplus_jM_j)$, then
$\Hom_A(V,M)\oplus L$ is an $n$-BB-tilting
$\End_A(V)$-module.\label{bbtilting}
\end{Prop}

{\it Proof.} The proof of (1) is similar to the one of Proposition
\ref{apr}. We leave it to the reader.

(2) We shall show that $L$ is isomorphic to
$\tau^{-1}\Omega_{\Lambda}^{-n+1}(S)$ with $S$ =
$\tau\Omega_{\Lambda}^{n-1}(L)$ being a simple $\Lambda$-module. It
is easy to see that $D(S)$ = Tr$\Omega_{\Lambda}^{n-1}(L)$ is a
simple right $\Lambda$-module. In fact, it is isomorphic to the top
of the indecomposable right $\Lambda$-module $\Hom_A(X,V)$, and is
not injective since $X\not\in\add(\bigoplus_j M_j)$.  Further, it
follows from $X\not\in\add(\bigoplus_i M_i)$ that we have an exact
sequence
$$\begin{array}{rl} 0\lra \Hom_A(Y,V)\lra\Hom_A(M_1,V) & \lra  \Hom_A(M_2,V)
\lra \cdots\lra\Hom_A(M_n,V)\\ & \\ \lra \Hom_A(X,V)   & \lra
\mbox{Tr}\Omega_{\Lambda}^{n-1}(L)=D(S)\lra 0.
\end{array}$$
If we apply $\Hom_{\Lambda}(-,\Lambda)$ to this sequence, we can see
that $\Ext^i_{\Lambda}(D(S),\Lambda)=0$ for all $0\le i\le n-1$.
This is just the condition (a) in the definition of an
$n$-BB-tilting module. To see that $\Ext^i_{\Lambda}(S,S)=0$ for all
$1\le i\le n$, we show that $\Ext^i_{\Lambda^{op}}(D(S),D(S))=0$ for
all $1\le i\le n$. This means that the projective cover
$\Hom_A(X,V)$ of the right $\Lambda$-module $D(S)$ is not a direct
summand of $\Hom_A(M_i,V)$ for all $1\le i\le n$. However, this
follows from the assumption that $X\not\in \add(\bigoplus_{j=1}^n
M_j)$. Thus the condition (b) of an $n$-BB-tilting module is
fulfilled. $\square$

\medskip
{\it Remarks.} (1) One can see that a non-injective simple
$A$-module $S$ defines an $n$-BB-tilting module if and only if
($a'$) proj.dim$_A(\tau^{-1}\Omega^{-n+1}(S))\le n$, ($b'$)
$\Ext^i_A(S,S)=0$ for all $1\le i\le n$ and ($c'$)
$\Ext^i_A(D(A),S)=0$ for all $1\le i\le n-1$. Note tat if a simple
module $S$ defines an $n$-BB-tilting module then the injective
dimension of $S$ is $n$ if and only if
$\Hom_A(\tau^{-1}\Omega^{-n+1}(S),A)=0$.

(2) With the same method as in Proposition \ref{bbtilting}, we can
prove the following fact:

Let $\cal C$ be a maximal $(n-1)$-orthogonal subcategory of $A$-mod
with $A$ a finite-dimensional algebra over a field $(n\ge 1)$.
Suppose $X$ and $Y$ are two indecomposable $A$-modules in $\cal C$
such that the sequence
$$0\lra X\stackrel{f}{\lra} M_n\stackrel{t_n}{\lra}M_{n-1}\lra\cdots
\lra M_2\stackrel{t_2}{\lra} M_1\stackrel{g}{\lra} Y\lra 0$$ is an
$n$-almost split sequence in $\cal C$. We define
$M=\bigoplus_{i=1}^n M_i$, $V = M\oplus X$, and $U= X\oplus M\oplus
Y$. Let $\Sigma$ be the endomorphism algebra of $U$. If $X\not\in
\add(M\oplus Y)$, then $T:= \Hom_A(V,U)\oplus S^X$ is an
(n+1)-BB-tilting right $\Sigma$-module, where $S^X$ is the top of
the right $\Sigma$-module $\Hom_A(X,U)$. If we define $\Delta
=\End(T_{\Sigma})$, then
$\Hom_{\Sigma}(\Hom_A(V,U)_{\Sigma},T_{\Sigma})\oplus
\Hom_{\Sigma}(\Hom_A(Y,U)_{\Sigma}, T_{\Sigma})$ is an
$(n+1)$-APR-tiling $\Delta$-module, that is, it is an
$(n+1)$-BB-tilting $\Delta$-module defined by the projective simple
$\Delta$-module $\Hom_{\Sigma}(S^X,T)$. Note that $\Delta$ is a
one-point extension of $\End_A(V)$ because
$\Hom_{\Sigma}(S^X,\Sigma)=0$.

\section{Auslander-Reiten triangles and derived equivalences
\label{sectriangle}}

By Corollary \ref{ars}, one can get a derived equivalence from an
Auslander-Reiten sequence. An analogue of an Auslander-Reiten
sequence in a triangulated category is the notion of
Auslander-Reiten triangle. Thus, a natural question rises: is it
possible to get a derived equivalence from an Auslander-Reiten
triangle in a triangulated category? In this section, we shall
discuss this question. First, let us briefly recall some basic
definitions concerning Auslander-Reiten triangles. For more details,
we refer the reader to \cite{HappelTri}.

Let $R$ be a commutative ring. Let $\mathcal{C}$ be a triangulated
$R$-category such that $\Hom_{\mathcal{C}}(X,Y)$ has finite length
as an $R$-module for all $X$ and $Y$ in $\mathcal{C}$. In this case,
we say that $\cal C$ is a Hom-finite triangulated $R$-category.
Suppose further that the category $\mathcal{C}$ is a Krull-Schmidt
category. A triangle
$X\stackrel{f}{\lra}M\stackrel{g}{\lra}Y\stackrel{w}{\lra}X[1]$ in
$\mathcal{C}$ is called an {\em Auslander-Reiten triangle} if

(AR1) $X$ and $Y$ are indecomposable;

(AR2) $w\ne 0$;

(AR3) if $t:U\lra Y$ is not a split epimorphism, then $tw=0$.

\medskip
Note that neither $f$ is a monomorphism nor $g$ is an epimorphism in
an Auslander-Reiten triangle. This is a difference of an
Auslander-Reiten triangle from an almost $\cal D$-split sequence.
Thus, an Auslander-Reiten triangle in a triangulated category may
not be an almost $\cal D$-split sequence. Also, an Auslander-Reiten
sequence in the module category of an Artin algebra in general may
not give us an Auslander-Reiten triangle in its derived module
category. For an Artin algebra, we even don't know whether its
stable module category has a triangulated structure except that the
Artin algebra is self-injective. In this case, an Auslander-Reiten
sequence can be extended to an Auslander-Reiten triangle in the
stable module category.

Recall that a morphism $f: U\lra V$ in a category $\mathcal{C}$ is
called a \emph{split monomorphism} if there is a morphism $g:V\lra
U$ in $\cal C$ such that $fg=id_U$; a \emph{split epimorphism} if
$gf=id_V$; and an \emph{irreducible} morphism if $f$ is neither a
split monomorphism nor a split epimorphism, and, for any
factorization $f=f_1f_2$ in ${\cal C}$, either $f_1$ is a split
monomorphism or $f_2$ is a split epimorphism.

Suppose
$X\stackrel{f}{\lra}M\stackrel{g}{\lra}Y\stackrel{w}{\lra}X[1]$
 is an Auslander-Reiten triangle in a triangulated category $\cal C$. Then we have the following
basic properties:

(1) $fg=0$ and $gw=0$. Moreover, both $f$ and $g$ are irreducible
morphisms.

(2) If $s: X\rightarrow U$ is not a split monomorphism,  then $s$
factors through $f$. Similarly, if $t: V\ra Y$ is not a split
epimorphism, then $t$ factors through $g$.

(3) Let $V$ be an indecomposable object in $\mathcal{C}$. Then $V$
is a direct summand of $M$ if and only if there is an irreducible
map from $V$ to $Y$ if and only if there is an irreducible map from
$X$ to $V$.

\medskip
We mention that in any triangulated category $\cal C$ the functors
$\Hom_{\cal C}(V,-)$ and $\Hom_{\cal C}(-,V)$ are co-homological
functors for each object $V\in {\cal C}$ (see \cite[Proposition 1.2,
p.4]{HappelTri}).

\medskip
The following is an expected result for Auslander-Reiten triangles.

\begin{Prop} Let $\mathcal{C}$ be a Hom-finite, Krull-Schmidt, triangulated
$R$-category. Suppose
$X\stackrel{f}{\lra}M\stackrel{g}{\lra}Y\stackrel{w}{\lra}X[1]$ is
an Auslander-Reiten triangle in $\mathcal{C}$ such that
$X[1]\not\in\add(M\oplus Y)$. If $N$ is an object in $\mathcal{C}$
such that none of $X, Y, X[1]$ and $ Y[-1]$ belongs to $\add(N)$,
then $\End_{\mathcal{C}}(N\oplus M\oplus X)$ and
$\End_{\mathcal{C}}(N\oplus M\oplus Y)$ are derived-equivalent via a
tilting module. In particular, $\End_{\mathcal{C}}(M\oplus X)$ and
$\End_{\mathcal{C}}(M\oplus Y)$ are derived-equivalent via a tilting
module. \label{artriangle}
\end{Prop}

{\it Proof.}  First, if $X$ is a direct summand of $M$, then there
is an  irreducible map from $X$ to $Y$. It follows from the property
(3) of an Auslander-Reiten triangle that $Y$ is a direct summand of
$M$. Similarly, if $Y$ is a direct summand of $M$, then so is $X$.
Thus, if $X$ or $Y$ is in $\add (M)$, then $\add(N\oplus M\oplus
X)=\add(N\oplus M\oplus Y)=\add (N\oplus M)$. In this case, both
$\End_{\mathcal{C}}(N\oplus M\oplus X)$ and
$\End_{\mathcal{C}}(N\oplus M\oplus Y)$ are Morita equivalent to
$\End_{\cal C}(N\oplus M)$, and therefore
$\End_{\mathcal{C}}(N\oplus M\oplus X)$ and
$\End_{\mathcal{C}}(N\oplus M\oplus Y)$ are derived-equivalent. Now,
we assume that neither $X$ nor $Y$ is in $\add(M)$. For simplicity,
we set $U:=N\oplus M$, $V:=U\oplus X$ and $W := U\oplus Y$. Denote
by $\Lambda$ the endomorphism ring of $V$. Since $X$ and $Y$ are not
in add$(U)$, we see that $f$ is a left add$(U)$-approximation of $X$
and $g$ is a right add$(U)$-approximation of $Y$. To see that the
condition (2) in Lemma \ref{3.1} is satisfied, we consider the exact
sequence
$$\cdots \ra \Hom_{\cal C}(V,M[-1])\stackrel{\delta}{\lra}
\Hom_{\cal C}(V,Y[-1])\ra \Hom_{\cal C}(V,X)\ra \Hom_{\cal
C}(V,M)\ra \Hom_{\cal C}(V,Y). $$ We have to show that the map
$\delta $ is surjective. By assumption, we have $Y[-1]\not\in
\add(N)$ and $Y[-1]\not\simeq X$ since $Y\not\simeq X[1]$. If
$Y[-1]\in\add(M)$, then there is an irreducible map from $X$ to
$Y[-1]$ by the property (3), and therefore there is an irreducible
map from $X[1]$ to $Y$. It follows that $X[1]$ is a direct summand
of $M$, which contradicts to our assumption that
$X[1]\not\in\add(M)$. This shows that $Y[-1]\not\in\add(M)$. Thus
any morphism from $V$ to $Y[-1]$ cannot be a split epimorphism. This
implies that the map $\delta$ is surjective by the property (2) of
an Auslander-Reiten triangle since the triangle $X[-1]\lra M[-1]\lra
Y[-1]\lra X $ is also an Auslander-Reiten triangle. Hence we have a
desired exact sequence
$$0\lra \Hom_{\cal C}(V,X)\lra \Hom_{\cal C}(V,M)\lra \Hom_{\cal
C}(V,Y). $$ Similarly, we have an exact sequence
$$0\lra \Hom_{\cal C}(Y,W)\lra \Hom_{\cal C}(M,W)\lra \Hom_{\cal C}(X,W). $$
Thus Proposition \ref{artriangle} follows from Lemma \ref{3.1} by
taking $n=1$. $\square$

\medskip
From Proposition \ref{artriangle} we get the following corollary.

\begin{Koro} Let $A$ be a self-injective Artin algebra.
Suppose $0\ra X\ra M\ra Y\ra 0$ is an Auslander-Reiten sequence such
that $\Omega^{-1}(X)\not\in\add(M\oplus Y)$. Then
$\underline{\End}_A(M\oplus X)$ and $\underline{\End}_A(M\oplus Y)$
are derived-equivalent, where $\underline{\End}_A(M)$ stands for the
quotient of $\End_A(M)$ by the ideal of those endomorphisms of $M$,
which factor through a projective $A$-module. \label{qftriangle}
\end{Koro}

{\it Proof.} If $A$ is a self-injective Artin algebra, then every
Auslander-Reiten sequence $0\ra X\ra M\ra Y\ra 0$ in $A$-mod can be
extended to an Auslander-Reiten triangle
$$ X\lra M\lra Y\lra \Omega_A^{-1}X$$ in the triangulated category
$A$-\underline{mod} which is equivalent to $\Db{A}/\Kb{A}$ (for
details, see \cite{HappelTri}). Thus Corollary \ref{qftriangle}
follows. $\square$

\medskip
Note that under the assumptions in Proposition \ref{artriangle} the
corresponding statement of Proposition \ref{apr} holds true for an
Auslander-Reiten triangle.

Finally, let us remark that Corollary \ref{qftriangle} may fail if
$A$ is not self-injective; for example, if we take $A$ to be the
path algebra (over a field $k$) of the quiver $2\lra 1\longleftarrow
3$, then there is an Auslander-Reiten sequence $$ 0\lra P(1)\lra
P(2)\oplus P(3)\lra I(1)\lra 0,$$ where $P(i)$ and $I(i)$ stand for
the projective and injective modules corresponding to the vertex
$i$, respectively. Clearly, this is a desired counterexample.

\section{An Example \label{sect7}}

In this section, we illustrate our results with an example.

\medskip
{\parindent=0pt\bf Example 1. } Let $k$ be a field, and let
$A=k[x,y]/(x^2,y^2)$. If $Y$ denotes the simple $A$-module, then
there is an Auslander-Reiten sequence $$0\lra X\lra N\oplus N\lra
Y\lra 0$$ in $A$-mod. Note that $X=\Omega_A^2(Y)$ and $N$ is the
radical of $A$. By Theorem \ref{thm1} or Corollary \ref{cor1}, the
two algebras $\End_A(N\oplus Y)$ and $\End_A(N\oplus X)$ are
derived-equivalent. Though the local diagram of the Auslander-Reiten
sequence is reflectively symmetric, the two algebras $\End_A(N\oplus
Y)$ and $\End_A(N\oplus X)$ are very different. This can be seen by
the following presentations of the two algebras given by quiver with
relations:

{\setlength{\unitlength}{0.7pt}
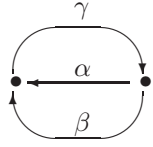
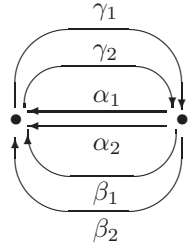
\begin{picture}(400,65)

\put(100,45){$\End_A(N\oplus Y)$}

\put(100,-40){$\bullet$}
\put(170,-40){$\bullet$}\put(137,-32){\oval(70,50)[t]}\put(137,-42){\oval(70,50)[b]}
\put(165,-37){\vector(-1,0){55}} \put(135,0){$\gamma$}
\put(135,-33){$\alpha$} \put(135,-63){$\beta$}
\put(100,-120){$\alpha\gamma=0=\beta\gamma.$}
\put(315,45){$\End_A(N\oplus X)$} \put(300,-40){$\bullet$}
\put(390,-40){$\bullet$}\put(348,-30){\oval(80,45)[t]}\put(350,-42){\oval(80,50)[b]}
\put(385,-40){\vector(-1,0){75}}\put(385,-32){\vector(-1,0){75}}
\put(348,-30){\oval(90,85)[t]}\put(348,-46){\oval(90,80)[b]}
\put(102,-48){\vector(0,1){5}}\put(172,-27){\vector(0,-1){5}}
\put(303,-51){\vector(0,1){5}}\put(310,-49){\vector(0,1){5}}
\put(393,-26){\vector(0,-1){5}}\put(388,-25){\vector(0,-1){5}}

\put(345, 20){$\gamma_1 $}\put(345,-2){$\gamma_2 $}
\put(345,-27){$\alpha_1 $}\put(345,-53){$\alpha_2 $}
\put(345,-80){$\beta_1 $}\put(345,-100){$\beta_2 $}
\put(300,-130){$\gamma_i\alpha_j=0=\gamma_i\beta_j, i\neq j,$}
\put(280,-145){$\gamma_1\beta_1=\gamma_2\beta_2, \quad
\gamma_1\alpha_1=\gamma_2\alpha_2,$}

\put(280,-160){$\alpha_1\gamma_2=\beta_1\gamma_1, \quad
\alpha_2\gamma_2=\beta_2\gamma_1.$}
\end{picture}}

\vspace{4.5cm} Note that the algebra $\End_A(N\oplus Y)$ is a
$7$-dimensional algebra of global dimension $2$, while the algebra
$\End_A(N\oplus X)$ is a $19$-dimensional algebra of global
dimension 3. Hence the two algebras are not stably equivalent of
Morita type since global dimension is invariant under stable
equivalences of Morita type (see \cite{Xi1}). A calculation shows
that the Cartan determinants of the both algebras equal $1$.

Recall that the Cartan matrix of an Artin algebra $A$ is defined as
follows: Let $S_1, \cdots, S_n$ be a complete list of non-isomorphic
simple $A$-modules, and let $P_i$ be a projective cover of $S_i$. We
denote the multiplicity of $S_j$ in $P_i$ as a composition factor by
$[P_i:S_j]$. The Catan matrix of $A$ is the $n\times n$ matrix
$([P_i:S_j])_{1\le i,j\le n}$,  and its determinant is called the
\emph{Cartan determinant} of $A$. It is well-known that the Cartan
determinant is invariant under derived equivalences.

\bigskip
{\bf Acknowledgements.} The authors thank I. Reiten and M.C.R.
Butler for comments, and Hongxing Chen at BNU for discussions on the
first version of the manuscript. Also, C.C.Xi thanks NSFC
(No.10731070) for partial support.

December 20, 2007; revised July 25, 2008.

\bigskip
Current address of W. Hu: School of Mathematical Sciences, Peking
University, Beijing 100871, P.R.China; email: huwei@math.pku.edu.cn.

\end{document}